\documentclass[11pt]{article}
\usepackage{amssymb,amsmath,latexsym}
\usepackage{tablists}
\usepackage{hyperref}
\usepackage[latin1]{inputenc}
\usepackage[french,english]{babel}
\usepackage{graphicx}
\newtheorem{ex}{Example}
\newtheorem{qu}{Question}

\newtheorem{thh}{Theorem }
\newtheorem{co}{Corollary}
\newcommand{\p}{{{\bf Proof.\,\,}}}

\newcommand{\biindice}[3]%
{

\begin{array}[t]{c}
#1\\
{\scriptstyle #2}\\
{\scriptstyle #3}
\end{array}

}

\begin{document}
\title{\Large  \textbf{Essential Descent Spectrum Equality }}
\author{\normalsize Abdelaziz Tajmouati$^{1}$, Hamid Boua$^{2}$ \\\\
\normalsize Sidi Mohamed Ben Abdellah University\\
\normalsize Faculty of Sciences Dhar El Mahraz \\
\normalsize Laboratory of Mathematical Analysis and Applications\\
\normalsize Fez, Morocco\\
\normalsize Email: abdelaziz.tajmouati@usmba.ac.ma\,\,\,\, hamid12boua@yahoo.com\\
}

\date{}\maketitle
\begin{abstract}
A bounded operator $T$ in a Banach space $X$ is said to satisfy the essential descent spectrum equality, if the descent spectrum of $T$ as an operator on $X$ coincides with the essential descent spectrum of $T$.
In this note, we give  some conditions under which the equality $\sigma_{desc}(T) = \sigma^e_{desc}(T)$ holds for a single operator $T$.
\end{abstract}
\section{\large Introduction}
Throughout this paper, $X$ denotes a complex Banach space and $\mathcal{B}(X)$ denotes the Banach algebra of all bounded linear
operators on $X$. Let $T\in \mathcal{B}(X)$, we denote by $T^*$,  $R(T)$, $N(T)$, $\rho(T)$, $\sigma(T)$ and $\sigma_{su}(T)$ respectively the adjoint, the range, the kernel, the resolvent set, the spectrum and the surjectivity spectrum of $T$. An operator $T\in \mathcal{B}(X)$ is called semi-regular if $R(T)$ is closed and $N(T^n)\subseteq R(T)$ for all positive integer $n$. The operator $T\in \mathcal{B}(X)$ is said to have the single-valued extension property at $\lambda_0 \in \mathbb{C}$,
abbreviated $T$ has the SVEP at $\lambda_0$, if for every neighbourhood $\mathcal{U}$ of $\lambda_0$ the
only analytic function $f : \mathcal{U} \rightarrow X$ which satisfies the equation $(\lambda I - T)f(\lambda) = 0$
is the constant function $f \equiv 0$. For an arbitrary operator $T \in \mathcal{B}(X)$ let $\mathcal{S}(T)= \{\lambda \in \mathbb{C}   :  T \mbox{ does not have the SVEP at } \lambda \}$. Note that $\mathcal{S}(T)$ is open and is contained in the interior of the point spectrum $\sigma_p(T)$.
The operator $T$ is said to have the SVEP if $\mathcal{S}(T)$ is empty.

For $T \in \mathcal{B}(X)$, the local resolvent set $\rho_{T}(x)$ of $T$ at the point $x \in X$ is defined as the set of all $\lambda \in \mathbb{C}$
for which there exist an open neighborhood $\mathcal{U}_{\lambda}$ of $\lambda$ and an analytic function $f : \mathcal{U}_{\lambda}\rightarrow X$ such that
$(T - \mu )f(\mu) = x \mbox{ for all } \mu \in \mathcal{U}_{\lambda}$
The local spectrum $\sigma_{T}(x)$ of $T$ at $x$ is then defined as $\sigma_{T}(x) = \mathbb{C}\backslash \rho_{T}(x)$. The local analytic solutions
occurring in the definition of the local resolvent set will be unique for all $x \in X$ if and only if $T$ has SVEP.
For every subset $F$ of $\mathbb{C}$, we define the local spectral subspace of $T$ associated with $F$ by
$X_T (F)= \{x \in X : \sigma_{T}(x) \subseteq F\}$.
Evidently, $X_T (F)$ is a hyperinvariant subspace of $T$, but not always closed.
A bounded linear operator $T\in \mathcal{B}(X)$ on a complex Banach
space $X$ is decomposable if every open cover $\mathbb{C} = U \cup V$ of the complex plane
$\mathbb{C}$ by two open sets $U$ and $V$ effects a splitting of the spectrum $\sigma(T)$ and of the
space $X$, in the sense that there exist $T$-invariant closed linear subspaces $Y$ and
$Z$ of $X$ for which $\sigma(T_{|Y}) \subseteq U$, $\sigma(T_{|Z}) \subseteq V$, and $X = Y + Z$.\\
A bounded linear operator $T$ on a complex Banach space $X$ has the decomposition property $(\delta)$ if
$X = \mathcal{X}_{T}(\overline{U}) + \mathcal{X}_{T}(\overline{V})$ for every open cover ${U, V}$ of $\mathbb{C}$. Where $\mathcal{X}_{T}(F)$ is the vector space of all elements $x\in X$ for which there exists an analytic function $f : \mathbb{C}\backslash F \rightarrow X$ such that $(T - \mu)f(\mu) = x$, for $\mu \in \mathbb{C}\backslash F$.

The descent of $T$ is $d(T)=\min\{q :  R(T^q)=R(T^{q+1})\}$, if no such $q$ exists, we let $d(T)=\infty$ \cite{Aie}, \cite{lay} and \cite{Mul}. The  descend spectrum is defined by: $\sigma_{desc}(T)=\{\lambda \in \mathbb{C}  :  d(\lambda -T)=\infty \}$.

Let $T\in \mathcal{B}(X)$, and consider the decreasing sequence
$c_n(T) := \dim(R(T^n)/R(T^{n+1}))$, $n\in \mathbb{N}$, see \cite{Gr}. Following M. Mbekhta and
M. M\"{u}ller \cite{Mb}, we shall say that $T$ has finite essential descent if $d_e(T) :=
\inf \{ n \in \mathbb{N} : c_n(T) < \infty \}$, where the infimum over the empty set is taken to
be infinite, is finite. This class of operators contains also every operator T such that the descent,
$d(T) = \inf \{ n \in \mathbb{N}  : c_n(T) = 0 \}$, is finite. In general $\sigma^e_{desc}(T)\subseteq \sigma_{desc}(T)$, and this inequality is strice. Indeed, consider the unilateral right shift operator $T$. According \cite{BKO}, $\sigma_{desc}(T)$ contains the closed unit disk, and  $\sigma^e_{desc}(T)$ is contained in the unit circle.

In \rm\cite{Bel} Olfa Bel Hadj Fredj has proved that, if $T \in \mathcal{B}(X)$ with a spectrum $\sigma(T)$ at most countable, then $\sigma_{desc}(T)=\sigma^e_{desc}(T)$. It is easy to construct an operator $T$ satisfying the essential descent spectrum equality such that $\sigma(T)$ is uncountable. For example, let $H$ be the Hilbert space $\ell^2(\mathbb{N})$ provided by the canonical basis $\{e_1,e_2,...\}$. Let $T\in \mathcal{B}(X)$ be
defined as $T(x_1,x_2,...) = (\frac{x_1}{2}, 0, x_2,x_3,...)$, $(x_n)_n\in \ell^2(\mathbb{N})$. From \cite{AI}, we have $\sigma_{su}(T) = \Gamma \cup \{\frac{1}{2}\}$, where $\Gamma$ denote the unit circle. Then $\sigma_{su}(T)$ of empty interior. According \cite{Bel}, $\sigma_{desc}(T)\backslash \sigma^e_{desc}(T)$ is an open. Since $\sigma_{desc}(T)\backslash \sigma^e_{desc}(T) \subseteq\sigma_{su}(T)$, therefore $\sigma_{desc}(T)= \sigma^e_{desc}(T)$. Motivated by the previous Example, our goal is to study the following question:
\begin{qu}
Let $T \in \mathcal{B}(X)$. If $\sigma(T)$ is uncountable,  under which condition on $T$ does
$\sigma_{desc}(T)=\sigma^e_{desc}(T)$ ?
\end{qu}

\section{Main results}
We start by the following results:
\begin{thh}\rm\cite{Bel}\label{aaa}
Let $T \in \mathcal{B}(X)$ be an operator for which $d_e(T)$ is finite.
Then there exists $\delta > 0$ such that for $0 < |\lambda| < \delta$ and $p := p(T)$, we have the following assertions:
\begin{enumerate}
\item $T - \lambda$ is semi regular;
\item $\dim N(T - \lambda)^n = n \dim N(T^{p+1})/\!\raisebox{-.65ex}{\ensuremath{N(T^{p})}}$ for all $n \in \mathbb{N}$;
\item codim$R(T -\lambda )^n = n\dim R(T^p)/\!\raisebox{-.65ex}{\ensuremath{R(T^{p+1}}}$ for all $n \in \mathbb{N}$.
\end{enumerate}
\end{thh}

\begin{co}\rm\cite{Bel}\label{a}
Let $T\in \mathcal{B}(X)$ be an operator of finite descent $d$. Then
there exists $\delta > 0$ such that the following assertions hold for $0 < |\lambda| < \delta$:
\begin{enumerate}
\item $T - \lambda$ is onto;
\item $\dim N(T - \lambda) = \dim N(T^{d+1})/\!\raisebox{-.65ex}{\ensuremath{N(T^{d})}}$.
\end{enumerate}
\end{co}

We have the following theorem.
\begin{thh}\label{aa}
Let $T \in \mathcal{B}(X)$, then
\begin{center}
$\sigma_{desc}(T)= \sigma^e_{desc}(T)\cup \overline{\mathcal{S}(T^*)}$
\end{center}
In particular, if $T^*$ has SVEP, then
\begin{center}
$\sigma_{desc}(T)=\sigma^e_{desc}(T) $
\end{center}
\end{thh}
\p If $\lambda \notin \sigma_{desc}(T)$, then  $\lambda \notin \sigma^e_{desc}(T)$. From \cite[Theorem 3.8]{Aie}, $T^*$ has the SVEP at $\lambda$, which
establish $\sigma^e_{desc}(T)\cup \mathcal{S}(T^*)\subseteq \sigma_{desc}(T)$, then it follows that $\sigma^e_{desc}(T)\cup \overline{\mathcal{S}(T^*)}= \overline{\sigma^e_{desc}(T)\cup \mathcal{S}(T^*)}\subseteq \overline{\sigma_{desc}(T)}=\sigma_{desc}(T)$.
For the other inclusion, let $\lambda$ be a complex number such that $T-\lambda$ has finite essential descent and $\lambda \notin \overline{\mathcal{S}(T^*)}$. According to theorem \ref{a}, there is $\delta > 0$, such that for $0 < |\lambda - \mu| < \delta$ and $p \in \mathbb{N}$, the operator $T - \mu$ is semi regular and codim  $R(T -\mu ) = \dim R(T-\lambda)^p/\!\raisebox{-.65ex}{\ensuremath{R(T-\lambda)^{p+1}}}$. Let $D^*(\lambda, \delta)=\{\mu : 0 < |\lambda - \mu| < \delta\}$. But since $\lambda \notin \overline{\mathcal{S}(T^*)}$, $D^*(\lambda, \delta)\backslash \mathcal{S}(T^*)$ is non-empty, and hence it follows that there exists $\mu_0 \in D^*(\lambda, \delta)$ such that $T-\mu_0$ is semi regular and $T^*$ has the SVEP at $\mu_0$. From \cite[Theorem 3.17]{Aie}, $T - \mu_0$ has finite descent. Also by Corollary \ref{aa}, there exists $\mu_1 \in D^*(\lambda, \delta)$ such that $T-\mu_1$ is surjective. Therefore codim$R(T -\mu )= \dim R(T-\lambda)^p/\!\raisebox{-.65ex}{\ensuremath{R(T-\lambda)^{p+1}}}=0$. It follows that $R(T-\lambda)^p=R(T-\lambda)^{p+1}$, which forces that $\lambda \notin \sigma_{desc}(T)$.

\begin{ex}

Suppose that $T$ is an unilateral weighted right shift
on $\ell^p(\mathbb{N})$, $1\leq p < \infty$, with weight sequence $(\omega_n)_{n\in \mathbb{N}}$. If
$c(T) = \lim_{n\rightarrow \infty} \inf(\omega_1...\omega_n)^{\frac{1}{n}} = 0$. By \rm\cite[Theorem 2.88]{Aie}, $T^*$ has SVEP, and from Theorem \ref{aa}, we have
$\sigma_{desc}(T)=\sigma^e_{desc}(T)$.
\end{ex}

A mapping $T : \mathcal{A} \rightarrow \mathcal{A}$ on a commutative complex Banach algebra $\mathcal{A}$ is said to
be a multiplier if
\begin{center}
$u(Tv) = (Tu)v \mbox{ for all } u, v \in \mathcal{A}.$
\end{center}
Any element $a \in \mathcal{A}$ provides an example, since, if $L_a : \mathcal{A} \rightarrow \mathcal{A}$ denotes the
mapping given by $L_a(u) := au$ for all $u \in \mathcal{A}$, then the multiplication operator La
is clearly a multiplier on $\mathcal{A}$. The set of all multipliers of $\mathcal{A}$ is denoted by $M(\mathcal{A})$. We recall that an algebra $\mathcal{A}$ is said to be semi-prime if $\{0\}$ is the only
two-sided ideal $J$ for which $J^2 = {0}$.
\begin{co}
Let $T \in M(\mathcal{A})$ be a multiplier on a semi-prime, regular and commutative Banach algebra $\mathcal{A}$
then
\begin{center}
$\sigma_{desc}(T)=\sigma^e_{desc}(T) $
\end{center}
\end{co}
\p
From \cite[Corollary 6.52]{Aie}, then $T^*$ has SVEP. Therefore by theorem \ref{aa}, $\sigma_{desc}(T)=\sigma^e_{desc}(T) $.

$T \in \mathcal{B}(X)$ is said to be supercyclic if it is only required that , for
some $x \in X$, the homogeneous orbit $\mathbb{C}.O(x,T)= \{\lambda T^n(x)  : n \in \mathbb{N} \}$ be dense in $X$.
\begin{co}
Let $T \in \mathcal{B}(X)$ be supercyclic operator. Then
\begin{center}
$\sigma_{desc}(T)=\sigma^e_{desc}(T)$
\end{center}
\end{co}

\p
If  $T \in \mathcal{B}(X)$ is supercyclic, according to  \cite[Proposition 1.26]{DYN}, either $\sigma_p(T^*) = \emptyset $ or $\sigma_p(T^*) = \{\lambda \}$ for some $\lambda \neq 0$, hence int$(\sigma_p(T^*))=\emptyset $, so $S(T^*)=\emptyset$, from Theorem \ref{aa}, we have $\sigma_{desc}(T)=\sigma^e_{desc}(T)$.

\begin{co}
Let $T \in \mathcal{B}(X)$. If $T$ satisfies one of the following properties:
\begin{enumerate}
\item $T$ is decomposable,
\item  $T$ has property $(\delta)$,
\end{enumerate}
then
\begin{center}
$\sigma_{desc}(T)=\sigma^e_{desc}(T) $
\end{center}
\end{co}
\p
 If $T$ is decomposable or $T$ has property $(\delta)$. We know from \cite[Theorem 1.2.7]{Lar} and \cite[Theorem 2.5.5]{Lar} that $T^*$ has SVEP, by
Theorem \ref{aa}, we have $\sigma_{desc}(T)=\sigma^e_{desc}(T) $.

\end{document}